\documentclass{amsart}

\usepackage{amssymb}
\theoremstyle{plain}

\newtheorem{thm}{Theorem}

\newtheorem{lem}[thm]{Lemma}

\theoremstyle{definition}
\newtheorem{df}[thm]{Definition}
\newtheorem{rem}[thm]{Remark}

\author{Michael Robinson} 

\title{Instability of a parabolic equation with a quadratic nonlinearity}

\begin{document}
\begin{abstract}
A nonlinear parabolic differential equation with a quadratic
nonlinearity is presented which has at least one equilibrium.  The
linearization about this equilibrium is asymptotically stable, but by
using a technique inspired by H. Fujita, we show that the equilibrium
is unstable in the nonlinear setting.  The perturbations used have the
property that they are small in every $L^p$ norm, yet they result in
solutions which fail to be global.
\end{abstract}


\maketitle

\section{Introduction}
This article demonstrates that in infinite-dimensional settings,
stability of the linearization about an equilibrium of a dynamical
system is not sufficient to ensure that the equilibrium is stable.
This is in stark contrast to the situation in finite-dimensional
settings, where stability of the linearized system implies stability
of the equilibrium.  (See \cite{BoyceDiPrima}, for instance.)  A
crucial point is that the system exhibited has a spectrum which
includes zero, so stability is possible (as in the unforced heat
equation), though not guaranteed.

We study classical solutions to the Cauchy problem
\begin{equation}
\label{pde}
\begin{cases}
\frac{\partial}{\partial t} u(t,x) = \frac{\partial^2}{\partial x^2} u(t,x) - 2 f(x) u(t,x) -
u^2(t,x)\\
u(0,x)=h(x) \in C^\infty(\mathbb{R})\\
t>0,x\in \mathbb{R},\\
\end{cases}
\end{equation}
where $f \in C_0^\infty(\mathbb{R})$ is a positive function with two
bounded derivatives.  (By $C_0^\infty$, we mean the space of smooth
functions which decay to zero.)  Since the linear portion of the right
side of \eqref{pde} is a sectorial operator, we can use it to define a
nonlinear semigroup. \cite{Henry} \cite{ZeidlerIIA} The standard
regularity theory for parabolic equations turns \eqref{pde} into a
smooth dynamical system, the behavior of which is largely controlled
by its equilibria.  This problem evidently has at least one
equilibrium, namely $u(t,x) = 0$ for all $t,x$.  Depending on the
exact choice of $f$, there may be other equilibria, however they will
not concern us here.  The linearized form of \eqref{pde} about this
equilibrium is evidently
\begin{equation}
\label{linearized_pde}
\frac{\partial}{\partial t} u(t,x) = \frac{\partial^2}{\partial x^2} u(t,x) - 2 f(x) u(t,x).
\end{equation}
The zero function is asymptotically stable for the linearized problem,
by a standard comparison principle argument. \cite{Evans} However,
using a technique pioneered by Fujita in \cite{Fujita}, we will show
that this equilibrium is not stable in the nonlinear problem, even if
the initial condition has small $p$-norm for every $1 \le p \le
\infty$.  Fujita showed that if $f \equiv 0$, then the zero function
is an unstable equilibrium of \eqref{pde}.  The cause of the
instability in \eqref{pde} is the decay of $f$, for if
$f=\text{const}>0$, then the comparison principle shows that the zero
function is stable.  We extend Fujita's result, so that roughly
speaking, since $f \to 0$ away from the origin, the system is less
stable to perturbations away from the origin.  Another indication that
there may be instability lurking (though not conclusive proof) is that
the decay of $f$ means that the spectrum of the linearized operator on
the right side of \eqref{linearized_pde} includes
zero. \cite{Mazya_2007}

\section{Motivation}
The problem \eqref{pde} describes a reaction-diffusion equation
\cite{FiedlerScheel}, or a diffusive logistic population model with a
spatially-varying carrying capacity.  The choice of $f$ positive means
that the equilibrium $u \equiv 0$ describes a population saturated at
its carrying capacity.  Without the diffusion term, this situation is
well known to be stable.  The decay condition on $f$ means that the
carrying capacity diminishes away from the origin.

The spatial inhomogeneity of $f$ makes the analysis of \eqref{pde}
much more complicated than that of typical reaction-diffusion
equations.  The existence of additional equilibria for \eqref{pde}
is a fairly difficult problem, which depends delicately on $f$.  (See
\cite{Brezis_1984} for a proof of existence of equilibria in a
related setting.)

\section{Instability of the equilibrium}
\label{instability_sec}

Given an $\epsilon>0$, we will construct an initial condition $h \in
C^\infty(\mathbb{R})$ for the problem \eqref{pde}, with $\|h\|_p <
\epsilon$ for each $1\le p \le \infty$, such that $\|u(t)\|_\infty \to
\infty$ as $t \to T < \infty$.  In particular, this implies that
$u\equiv 0$ is not a stable equilibrium of \eqref{pde}, at least
insofar as classical solutions are concerned.  We employ a technique
of Fujita, which provides sufficient conditions for equations like
\eqref{pde} to blow up. \cite{Fujita} (Additionally, \cite{Evans}
contains a more elementary discussion of the technique with a similar
construction.)  Our choice for $h$ can be thought of as a sequence of
progressively shifted gaussians, and we will demonstrate that though
each has smaller $p$-norm than the previous, the solution started at
$h$ still blows up.

\subsection{The technique of Fujita}

The technique of Fujita examines the blow-up behavior of nonlinear
parabolic equations by treating them as ordinary differential equations
on a Hilbert space.  Suppose $u(t)$ solves 
\begin{equation}
\label{fujita_eqn}
\frac{\partial u(t)}{\partial t} = L u(t) + N(u(t),t),
\end{equation}
where $L$ is a linear operator not involving $t$, and $N$ may be
nonlinear and may depend on $t$.  Suppose that $v(t)$ solves 
\begin{equation}
\label{fujita_adj_eqn}
\frac{\partial v(t)}{\partial t} = - L^* v(t),
\end{equation}
where $L^*$ is the adjoint of $L$.  Let $J(t) = \left<v(t),u(t)\right>$.  We
observe that if $|J(t)| \to \infty$ then either $\|v(t)\|$ or
$\|u(t)\|$ also does.  So if $v(t)$ does not blow up, then we can show
that $\|u(t)\|$ blows up, and perhaps more is true.  If we
differentiate $J(t)$, we obtain the identity
\begin{eqnarray*}
\frac{d}{dt}J(t) &=& \frac{d}{dt} \left< v(t), u(t) \right>\\
&=& \left<\frac{dv}{dt},u(t)\right>+\left<v(t),\frac{du}{dt}\right>\\
&=&\left<-L^*v(t),u(t)\right>+\left<v(t),Lu(t)+N(u(t),t)\right>\\
&=&\left<v(t),N(u(t),t)\right>, 
\end{eqnarray*}
where there is typically a technical justification required for the
second equality.  It is often possible to find a bound for
$\left<v(t),N(u(t),t)\right>$ in terms of $J(t)$.  So then the method
provides a fence (in the sense of \cite{HubbardWest}) for $J(t)$,
which we can solve to give a bound on $|J(t)|$.  As a result, the
blow-up behavior of $u(t)$ is controlled by the solution of an {\it
ordinary} differential equation (for $J(t)$) and a {\it linear}
parabolic equation (for $v(t)$), both of which are much easier to
examine than the original nonlinear parabolic equation.

\subsection{Instability in $L^p$ for $1 \le p \le \infty$}

We begin our application of the method of Fujita by working with
$L=\frac{\partial^2}{\partial x^2}-2f$ and $N(u)=- u^2$ in
\eqref{fujita_eqn}.  Since \eqref{fujita_adj_eqn} is then not
well-posed for all $t$, we must be a little more careful than the
method initially suggests.  For this reason, we consider a family of
solutions $v_\epsilon$ to \eqref{fujita_adj_eqn} that have slightly
extended domains of definition.  It will also be important, for technical
reasons, to enforce the assumption that the first and second
derivatives of $f$ are bounded.  

\begin{df}
\label{v_eps_2_df}
Suppose $w=w(t,x)$ solves
\begin{equation}
\label{w_eqn}
\begin{cases}
\frac{\partial w}{\partial t} = \frac{\partial^2 w}{\partial x^2} - 2
f(x) w(t,x)\\
w(0,x) = w_0(x) \ge 0.
\end{cases}
\end{equation}
Define $v_\epsilon(s,x) = w(t-s+\epsilon,x)$ for fixed $t>0$ and
$s<t+\epsilon$.  Notice that by the comparison principle,
$v_\epsilon(s,x) \ge 0$.
\end{df}

\begin{lem}
\label{w_endcond_lem}
Suppose that $w$ solves \eqref{w_eqn}.  Then $w, \frac{\partial
  w}{\partial x} \in C_0(\mathbb{R})$.
\begin{proof}
The standard existence and regularity theorems for linear parabolic
equations (see \cite{ZeidlerIIA}, for example) give that
$w,\frac{\partial w}{\partial x},\frac{\partial^2 w}{\partial x^2} \in
L^2(\mathbb{R})$ and that $w \in C^2(\mathbb{R})$.  The comparison
principle, applied to $\frac{\partial}{\partial t}\frac{\partial
w}{\partial x}$ and $\frac{\partial}{\partial t}\frac{\partial^2
w}{\partial x^2}$ gives that the first and second derivatives of $w$
are bounded for each fixed $t$.  (This uses our assumption that $f$
has two bounded derivatives.)

The lemma follows from a more general result: if $g\in C^1 \cap L^p
(\mathbb{R})$ for $1 \le p < \infty$ and $g' \in
L^\infty(\mathbb{R})$, then $g \in C_0(\mathbb{R})$.  To show this, we
suppose the contrary, that $\lim_{x \to \infty} g(x) \ne 0$ (and
possibly doesn't exist).  By definition, this implies that there is
an $\epsilon >0$ such that for all $x>0$, there is a $y$ satisfying
$y>x$ and $|g(y)|>\epsilon$.  Let $S=\{y|\,|g(y)|>\epsilon\}$, which
is a union of open intervals, is of finite measure, and has $\sup S =
\infty$.  Let $T = \{y | \, |g(y)|>\epsilon/2\}$.  Note that $T$
contains $S$, but since $g'$ is bounded, for each $x \in S$, there is
a neighborhood of $x$ contained in $T$ of measure at least
$\epsilon/\|g'\|_\infty$.  Hence, since $\sup T = \sup S = \infty$,
$T$ cannot be of finite measure, which contradicts the fact that $g
\in L^p(\mathbb{R})$ with $1\le p < \infty$.
\end{proof}
\end{lem}

\begin{lem}
\label{fence_2_lem}
Suppose $u:[0,T)\times \mathbb{R} \to \mathbb{R}$ is a classical
  solution to \eqref{pde} with $u \le 0$ and $u(t)\in
  L^\infty(\mathbb{R})$ for each $t \in [0,T)$.  Then
\begin{equation}
\label{fence_2_eqn}
-\int w(t,x) h(x) dx \le \left( \int_0^t \frac{1}{\|w(s)\|_1}
 ds \right)^{-1},
\end{equation}
where $w$ is defined as in Definition \ref{v_eps_2_df}.

\begin{proof}
Define 
\begin{equation}
J_\epsilon(s)=\int{v_\epsilon(s,x)u(s,x) dx}.
\end{equation}

First of all, we observe that since $u\in L^\infty(\mathbb{R})$,
  $v_\epsilon(s,\cdot)u(s,\cdot)$ is in $L^1(\mathbb{R})$ for each
$s<t$.

Now suppose we have a sequence $\{m_n\}$ of compactly supported
smooth functions with the following properties: \cite{LeeSmooth}
\begin{itemize}
\item $m_n \in C^\infty(\mathbb{R})$,
\item $m_n(x) \ge 0$ for all $x$, \item $\text{supp}(m_n)$ is contained in the interval $(-n-1,n+1)$, and
\item $m_n(x)=1$ for $|x| \le n$.
\end{itemize}
Then it follows that 
\begin{equation*}
J_\epsilon(s)=\lim_{n\rightarrow \infty} \int{v_\epsilon(s,x)u(s,x)m_n(x) dx}.
\end{equation*}

Now 
\begin{eqnarray*}
\frac{d}{ds} J_\epsilon(s) &=& \frac{d}{ds} \lim_{n\rightarrow \infty}
\int{v_\epsilon(s,x)u(s,x)m_n(x) dx} \\
&=& \lim_{h\rightarrow 0} \lim_{n\rightarrow \infty} \frac{1}{h}\int
\left(v_\epsilon(s+h,x)u(s+h,x)-v_\epsilon(s,x)u(s,x)\right)m_n(x) dx.\\
\end{eqnarray*}
We'd like to exchange limits using uniform convergence.  To do this we
show that 
\begin{equation}
\label{big_lim}
\lim_{n\rightarrow \infty} \lim_{h\rightarrow 0} \frac{1}{h}\int
\left(v_\epsilon(s+h,x)u(s+h,x)-v_\epsilon(s,x)u(s,x)\right)m_n(x) dx
\end{equation}
exists and the inner limit is uniform.  We show both together by a
little computation, using uniform convergence and LDCT:

\begin{eqnarray*}
&&\lim_{n\rightarrow \infty} \lim_{h\rightarrow 0} \frac{1}{h}\int
\left(v_\epsilon(s+h,x)u(s+h,x)-v_\epsilon(s,x)u(s,x)\right)m_n(x) dx\\
&=&\lim_{n\rightarrow \infty} \int
\left(\frac{\partial}{\partial s}v_\epsilon(s,x)u(s,x)+v_\epsilon(s,x)\frac{\partial}{\partial s}u(s,x)\right)m_n(x) dx \\
&=&\lim_{n\rightarrow \infty} \int
\left(-\frac{\partial^2}{\partial x^2}
v_\epsilon(s,x)+2f(x)v_\epsilon(s,x)\right) u(s,x) m_n(x) +\\
&&v_\epsilon(s,x)\left(\frac{\partial^2}{\partial x^2} u(s,x) -
u^2(s,x)- 2 f(x) u(x) \right)m_n(x) dx \\
&=&\lim_{n\rightarrow \infty} \int
-v_\epsilon(s,x)u^2(s,x)m_n(x) dx.
\end{eqnarray*}
Minkowski's inequality has that
\begin{equation*}
\left|\int v_\epsilon u m_n dx \right|\le \int v_\epsilon |u| m_n dx \le \left(\int v_\epsilon m_n dx\right)^{1/2}
\left(\int v_\epsilon u^2 m_n dx \right)^{1/2},
\end{equation*}
since $v_\epsilon, m_n\ge 0$.  This gives that

\begin{eqnarray*}
&&\int
-v_\epsilon(s,x)u^2(s,x)m_n(x) dx\\
&\le& - \frac{(\int v_\epsilon u m_n dx)^2 }{\int v_\epsilon m_n dx}\\
&\le& - \frac{\left(\int v_\epsilon u dx \right)^2}{\int v_\epsilon m_1 dx}, \\
\end{eqnarray*}
hence the inner limit of \eqref{big_lim} is uniform.  On the other
hand, 
\begin{equation*}
|v_\epsilon(s,x)u^2(s,x)m_n(x)| \le v_\epsilon(s,x) \|u(s)\|_\infty^2
 \in L^1(\mathbb{R})
\end{equation*}
so the double limit of \eqref{big_lim} exists by dominated
convergence.  Thus we have the fence
\begin{equation}
\label{jeps_2_eqn}
\frac{d J_\epsilon(s)}{ds} \le -\frac{(J_\epsilon(s))^2}{\|v_\epsilon(s)\|_1}.
\end{equation}

We solve the fence \eqref{jeps_2_eqn} to obtain (note $J_\epsilon \le 0$)
\begin{eqnarray*}
\frac{1}{\|v_\epsilon(s)\|_1} &\le& - \frac{d J_\epsilon(s)}{ds}
\frac{1}{(J_\epsilon(s))^2}\\
\int_0^t \frac{1}{\|v_\epsilon(s)\|_1} ds &\le& \frac{1}{J_\epsilon(t)}
- \frac{1}{J_\epsilon(0)}\\
\int_0^t \frac{1}{\|v_\epsilon(s)\|_1} ds &\le& -\frac{1}{J_\epsilon(0)}.
\end{eqnarray*}
Taking the limit as $\epsilon \to 0$ of both sides of the inequality
yields
\begin{equation*}
-\int w(t,x) h(x) dx \le \left( \int_0^t \frac{1}{\|w(t-s)\|_1} ds
 \right)^{-1}
 =  \left( \int_0^t \frac{1}{\|w(s)\|_1} ds \right)^{-1},
\end{equation*}
as desired.
\end{proof}
\end{lem}

\begin{rem}
\label{l_inf_stability}
Since we are interested in proving the instability of the zero
function in \eqref{pde}, consider $u(0,x) = h(x) = -
\epsilon$ for $\epsilon > 0$.  Then \eqref{fence_2_eqn} takes on the
simple form
\begin{equation}
\label{fence_2_improved_easy}
\epsilon \int_0^t \frac{\|w(t)\|_1}{\|w(s)\|_1} ds \le 1.
\end{equation}
So in particular, $\|u(t)\|_\infty$ blows up if there exists a $T>0$
such that $\epsilon \int_0^T \frac{\|w(T)\|_1}{\|w(s)\|_1} ds > 1.$

The stability of the zero function in \eqref{pde} depends
on the stability of the zero function in \eqref{w_eqn} -- the
linearized problem.  If the zero function in the linearized problem is
very strongly attractive, say $\|w(t)\|_1 \sim e^{-t}$, then
\begin{equation*}
\int_0^t \frac{e^{-t}}{e^{-s}} ds = (1-e^{-t}) < 1,
\end{equation*}
and so a small choice of $\epsilon<1$ does not cause blow-up via a
violation of \eqref{fence_2_improved_easy}.  On the other hand,
blow-up occurs if it is less attractive, say $\|w(t)\|_1 \sim
t^{-\alpha}$ for $\alpha \ge 0$.  Because then
\begin{equation*}
\int_0^t \frac{s^\alpha}{t^\alpha} ds = \frac{t}{\alpha + 1},
\end{equation*}
whence blow-up occurs before $t=\frac{\alpha + 1}{\epsilon}$.

In the particular case of $f(x)=0$ for all $x$, we note that $w$ is
simply a solution to the heat equation, which has
$\|w(t)\|_1=\|w_0\|_1$ for all $t$ (by direct computation
using the fundamental solution, say), so blow up occurs.  Thus
we can recover a special case of the original blow-up result of Fujita
in \cite{Fujita}.
\end{rem}

\begin{thm}
\label{fence_2_improved_violation}
Suppose a sufficiently small $\epsilon>0$ is given.  Then for a
certain choice of initial condition $h(x)$ with
$\|h\|_p < \epsilon$ for all $1\le p \le \infty$, there exists a
$T>0$ for which $\lim_{t\to T^-}\|u(t)\|_\infty = \infty$.
\begin{proof}
First, it suffices to choose $\|u(0)\|_1<\epsilon$ and $\|u(0)\|_\infty<\epsilon$, since
\begin{equation*}
\|u\|_p = \left(\int |u|^p dx\right)^{1/p} \le \|u\|_\infty^{(p-1)/p}
\|u\|_1^{1/p} < \epsilon.
\end{equation*}  
We assume, contrary to what is to be proven, that $\|u(t)\|_\infty$
does not blow up for any finite $t$.  In other words, assume that $u:
[0,\infty)\times \mathbb{R} \to \mathbb{R}$ is a classical solution to
\eqref{pde}, with $\|u(t)\|_\infty < \infty$ for all $t$.  We make
several definitions:
\begin{itemize}
\item Choose $0<\beta < \min\left\{\epsilon,\frac{\epsilon^4}{16 \pi^2}\right\}$.
\item Choose $\gamma>0$ small enough so that 
\begin{equation}
\label{vio_1}
\frac{\beta}{27 \gamma^2} = K,
\end{equation}
for some some arbitrary $K > 1$.
\item Since $0 \le f \in C_0^\infty(\mathbb{R})$, we can choose an
  $x_1$ such that
\begin{equation}
\label{vio_2}
f(x) \le \gamma \text{ when } x < x_1.
\end{equation}
\item Next, we choose $x_0<x_1$ so that
\begin{equation}
\label{vio_3}
\sqrt{t}\|f\|_\infty \left ( 1 - \text{erf} \left (
\frac{x_1-x_0}{2\sqrt{t}} \right ) \right ) < \gamma 
\end{equation}
for all $0 < t < \frac{1}{4 \gamma^2}$.  Notice that any choice less
than $x_0$ will also work.
\item Choose the initial condition for \eqref{pde} to be
\begin{equation}
\label{initial_choice_eqn}
u(0,x)=h(x)=- \beta e^{\beta^{3/2} (x-x_0)^2}.
\end{equation}
This choice of initial condition has $\|u(0)\|_\infty =
\beta<\epsilon$, $\|u(0)\|_1 = 2 \pi^{1/2} \beta^{1/4}<\epsilon$, and $\left \|
\frac{\partial^2 u(0)}{\partial x^2}\right \|_\infty = \mu = 2
\beta^{5/2}.$  (The value of $\mu$ will be important shortly.)
\item Finally, let $w_0(y)=\delta(y-x_0)$ (the Dirac
  $\delta$-distribution), and suppose that $w$ solves \eqref{w_eqn}.
  In other words, choose $w$ to be the fundamental solution to
  \eqref{w_eqn} concentrated at $x_0$.  Note that the maximum
  principle ensures both that $w(t,x) \ge 0$ for all $t>0$ and $x \in
  \mathbb{R}$ and that $\|w(t)\|_1 \le \|w(0)\|_1 = 1$ for all $t>0$.
  This allows us to rewrite \eqref{fence_2_eqn} as
\begin{equation}
\label{fence_2_improved}
- t \int w(t,x) h(x) dx \le 1.
\end{equation}
\end{itemize}

Now we estimate the integral in \eqref{fence_2_improved}.  Notice that
\begin{eqnarray*}
\frac{d}{dt}\int w(t,x) \left(-h(x)\right) dx &=& 
\int \left(\frac{\partial^2 w}{\partial x^2} - 2 f(x) w(t,x)\right)
\left(-h(x)\right) dx \\
&=&\int \left(-\frac{\partial^2 u}{\partial x^2} + 2 f(x) h(x) \right) w(t,x)dx, \\
\end{eqnarray*}
where Lemma \ref{w_endcond_lem} eliminates the
boundary terms.  Now suppose $z$ solves the heat equation with the
same initial condition as $w$, namely
\begin{equation}
\label{z_eqn}
\begin{cases}
\frac{\partial z}{\partial t} = \frac{\partial^2 z}{\partial x^2}\\
z(0,x) = w_0(x) = \delta(x-x_0).
\end{cases}
\end{equation}
The comparison principle estabilishes that $z(t,x) \ge w(t,x)$ for all
$t>0$ and $x \in \mathbb{R}$, since $f,w \ge 0$.  As a result, we have
that
\begin{eqnarray*}
\frac{d}{dt}\int w(t,x) \left(-h(x)\right) dx &\ge&
\int \left(-\left |\frac{\partial^2 u}{\partial x^2}\right | + 2 f(x) h(x) \right)
z(t,x)dx \\
&\ge&
- \mu - 2 \beta \int f(x) z(t,x) dx,\\
\end{eqnarray*}
where $\mu=\left\|\frac{\partial^2 u}{\partial x^2}(0)\right\|_\infty$ and
$\beta=\|u(0)\|_\infty$, which is an integrable equation.  As a result,
\begin{equation}
\label{w_1_bnd}
\int w(t,x) \left(-h(x)\right) dx \ge \beta -\mu t - 2 \beta \int_0^t
\int \int f(x) \frac{1}{\sqrt{4 \pi s}} 
e^{-\frac{(x-y)^2}{4s}}w_0(y) dy\, dx\, ds.
\end{equation}
On the other hand using our choice for $w_0$,
\begin{eqnarray*}
\int_0^t \int \int & f(x)& \frac{1}{\sqrt{4 \pi s}}
e^{-\frac{(x-y)^2}{4s}}w_0(y) dy\, dx\, ds = \int_0^t \int f(x) \frac{1}{\sqrt{4 \pi s}}
e^{-\frac{(x-x_0)^2}{4s}} dx\, ds\\
&\le& \int_0^t \frac{1}{\sqrt{4\pi s}} \left(\gamma
\int_{-\infty}^{x_1} e^{-\frac{(x-x_0)^2}{4s}} dx + \|f\|_\infty
\int_{x_1}^\infty e^{-\frac{(x-x_0)^2}{4s}} dx \right) ds\\
&\le& \frac{\gamma \sqrt{t}}{4} + \frac{1}{2}\|f\|_\infty \int_0^t 1 -
\text{erf}\left( \frac{x_1-x_0}{2\sqrt{s}}\right) ds\\
&\le& \frac{\gamma \sqrt{t}}{4} + \frac{1}{2}\|f\|_\infty \int_0^t 1 -
\text{erf}\left( \frac{x_1-x_0}{2\sqrt{t}}\right) ds\\
&\le& \frac{\gamma \sqrt{t}}{4} + \frac{1}{2}t \|f\|_\infty \left( 1 - \text{erf}
\left( \frac{x_1-x_0}{2\sqrt{t}} \right ) \right )\\
&\le& \frac{3 \gamma \sqrt{t}}{4} \le \gamma \sqrt{t},\\
\end{eqnarray*}
we have used \eqref{vio_2}, \eqref{vio_3}, and assumed that $0 < t <
\frac{1}{4 \gamma^2}$.  Then \eqref{fence_2_improved} becomes
\begin{equation*}
1 \ge t \int w(t,x) \left(-h(x)\right) dx \ge \beta t - \mu t^2 - 2
\beta \gamma t \sqrt{t} = -2 \beta^{5/2} t^2 - \frac{2 \beta^{3/2}
  t^{3/2}}{\sqrt{27 K}} + \beta t,
\end{equation*}
using our choices of $\mu$, $\gamma$, and initial condition.  Maple
reports that the maximum of $A(t)=-2 \beta^{5/2} t^2 - \frac{2 \beta^{3/2}
  t^{3/2}}{\sqrt{27 K}} + \beta t$ is unique, occurs at $0 <
t_0 < \frac{1}{4 \gamma^2}$, and has the asymptotic expansion
\begin{equation*}
A(t_0) \sim K - 18 K \sqrt{\beta} + 432 K^3 \beta + O(\beta^{3/2}).
\end{equation*}
Thus for all small enough $\epsilon > \beta$, we obtain a
contradiction to \eqref{fence_2_improved} since $K>1$.  Thus, for some
$T<t_0<\infty$, $\lim_{t\to T^-} \|u(t)\|_\infty = \infty$.
\end{proof}
\end{thm}

\section{Discussion}
Theorem \ref{fence_2_improved_violation} gives a fairly strong
instability result.  No matter how small an initial condition to
\eqref{pde} is chosen, even with all $p$-norms chosen small, solutions
can blow up so quickly that they fail to exist for all $t$.  This
precludes any kind of stability for classical solutions.  Like the
analogous result in Fujita's paper, the kind of initial conditions
which can be responsible for blow up are of the nicest kind imaginable
-- gaussians in either case!  

It must be understood that the argument in Theorem
\ref{fence_2_improved_violation} depends crucially on the decay of
$f$.  Without it, the lower bound on $\int w(t,x) (-h(x)) dx$
decreases too quickly.  Indeed, if $f=\text{const}>0$ and $h(x)>-f$,
then the comparison principle demonstrates that the zero function is
asymptotically stable.  On the other hand, any rate of decay for $f$
satisfies the hypotheses of Theorem \ref{fence_2_improved_violation},
and so will cause \eqref{pde} to exhibit instability.

Finally, although we have examined the case where the nonlinearity in
\eqref{pde} is due to $u^2$, there is no obstruction to extending the
analysis to any nonlinearity like $|u|^k$, with degree $k$ greater than
2.  A higher-degree nonlinearity would result in a somewhat different
form for \eqref{fence_2_eqn}, but this presents no further
difficulties to the argument.  Indeed, by analogy with Fujita's work,
higher-degree nonlinearities would result in significantly faster
blow-up.

\bibliography{instability_bib}
\bibliographystyle{hplain}

\end{document}